\documentclass[11pt]{amsart}
\usepackage{amsfonts,amssymb,amsthm,amsmath,amsxtra,amscd,verbatim,eucal}

\usepackage[all]{xy}
\usepackage[dvips]{graphics}

\title{Mather discrepancy and the arc spaces}
\author{Shihoko Ishii} 
\address{Graduate School of Mathematical Science,
University of Tokyo, Komaba, Meguro, Tokyo, Japan
\newline
e-mail :  shihoko@ms.u-tokyo.ac.jp}

\newcommand{\bC}{{\Bbb C}}
\newcommand{\bP}{{\Bbb P}}

\newcommand{\bN}{{\Bbb N}}
\newcommand{\bA}{{\Bbb A}}

\newcommand{\codim}{\operatorname{codim}}
\newcommand{\spec}{\operatorname{Spec}}

\newcommand{\sing}{\operatorname{Sing}}
\newcommand{\ord}{\operatorname{ord}}
\newcommand{\val}{\operatorname{val}}

\newcommand{\sT}{{\spec K[[t]]}}

\let \cedilla =\c
\renewcommand{\c}[0]{{\mathbb C}}  

\renewcommand{\o}[0]{{\mathcal O}} 

\newcommand{\n}[0]{{\mathbb N}}

\newcommand{\isom}{{\widetilde{\to}}}

\newcommand{\hx}{{\widehat X}}
\newcommand{\hk}{{\widehat K}}
\renewcommand{\a}{{\frak{a}}}

\newcommand{\lct}{\operatorname{lct}}
\newcommand{\mld}{\operatorname{mld}}
\renewcommand{\j}{{\mathcal J}}
\newcommand{\hlct}{\widehat{\operatorname{lct}}}
\newcommand{\hmld}{\widehat{\operatorname{mld}}}

\newcommand{\hKY}{\widehat K_{Y/X}}
\newcommand{\KY}{K_{Y/X}}
\newcommand{\im}{{\operatorname{Im}}}
\newcommand{\q}[0]{{\mathbb Q}}

\newcommand{\cont}{\operatorname{Cont}}

\def\to {\longrightarrow}

\newtheorem{thm}{Theorem}[section]

\newtheorem{lem}[thm]{Lemma}
\newtheorem{cor}[thm]{Corollary}
\newtheorem{prop}[thm]{Proposition}

\theoremstyle{definition}
\newtheorem{defn}[thm]{Definition}

\newtheorem{say}[thm]{}
\newtheorem{exmp}[thm]{Example}
\newtheorem{conj}[thm]{Conjecture}
 
\newtheorem{rem}[thm]{Remark}

\theoremstyle{remark}

\begin{document}
\maketitle
\footnote{partially supported by JSPS Grant-In-Aid (B) 22340004, (S)
19104001}

\begin{abstract}
  A goal of this paper is a characterization of  singularities according to 
  a new invariant, Mather discrepancy. 
 We also show some evidences convincing us that 
  Mather discrepancy is a reasonable invariant in a view point of birational geometry.

\end{abstract}

\section{Introduction}

\noindent
  Let $Y\to X$ be a resolution of the singularities of a variety $X$ over $\c$ that factors through the Nash blow-up $\hx\to X$.
  In \cite{DEI} we introduce the Mather discrepancy divisor $\hk_{Y/X}$ and showed its  relation  with the geometry of the arc space $X_\infty$ of $X$.
  Note that $\hk_{Y/X}$ coincides with the usual discrepancy divisor $K_{Y/X}$ 
  when $X$ is non-singular.
  For a pair $(X, \a)$ consisting of a $\q$-Gorestein variety $X$ and an ideal $\a\subset \o_X$, the log-canonical threshold $\lct(X, \a)$ and minimal log-discrepancy $\mld(X,\a)$ are defined by using the usual discrepancy divisor $K_{Y/X}$, where 
 $Y\to X$ is a log-resolution of  $(X, \a)$.
 In this paper, for a pair $(X, \a)$ consisting of an arbitrary variety $X$ and an ideal 
 $\a\subset \o_X$, we define  Mather log-canonical threshold $\hlct(X, \a)$ and Mather minimal log-discrepancy $\hmld(X,\a)$  by using the Mather discrepancy divisor $\hk_{Y/X}$, where 
 $Y\to X$ is a log-resolution of  $(X, \a)$ factoring through the Nash blow-up
 $\hx\to X$.
 These invariants are defined for an arbitrary variety $X$ without the condition of 
 $\q$-Gorenstein or even normality.
  On  these invariants, we obtain some formulas  which are generalizations of known formulas. 
  This shows that the Mather discrepancy is a reasonable invariant and we can expect that it will be useful in birational geomety. 
  The main result of this paper is about $\hmld(x; X, \o_X)$.
  
For a closed point $x\in X$ of a variety of dimension $n$, we have the following 
inequality:
$$\hmld(x; X, \o_X)\geq n.$$
We can characterize a singularity $(X,x)$  with $\hmld(x; X, \o_X)= n$.
Let $\lambda_m$ be defined by
$\dim\pi_m^{-1}(x)\cap \psi_m(X_\infty)=mn-\lambda_m$ and let $\lambda_m^0$ be 
defined by $\dim \pi_m^{-1}(x)\cap \psi_m(X_\infty\setminus (\sing X)_\infty)=mn-\lambda_m^0$,
where $\psi_m: X_\infty\to X_m$ and $\pi_m:X_m\to X$ are the truncation morphisms. Then it turns out that $\lambda_m$ and $\lambda_m^0$ are constant for $m\gg 0$ and they coincide with $\hmld(x; X, \o_X)- n$.
Our characterization of a singularity $(X,x)$  with $\hmld(x; X, \o_X)= n$ is  the following: 
\begin{thm}
\label{hmld=n}
For a singularity $(X,x)$ the following are equivalent:
\begin{enumerate}
\item[(i)] $\hmld(x; X, \o_X)= n$;
\item[(ii)] $\lambda_m=0$ for every $m\in \n$;
\item[(iii)] $\lambda^0_m=0$ for every $m\in \n$;
\item[(iv)] $\lambda^0_1=0$;
\item[(v)] The tangent cone of $(X,x)$ has a reduced irreducible component.
\end{enumerate}
\end{thm}

This paper is organized as follows: In the second section we put basic facts on Mather discrepancy and the arc spaces, which will be used in this paper.
In the third section the invariants $\hlct$ and $\hmld$ are defined by using Mather discrepancy and we prove 
the formulas on these invariants.
These formulas are natural generalization of the known formulas.
In the forth section we prove Theorem \ref{hmld=n} .

 We assume that every scheme in this paper is defined over $\c$.
  A variety is an irreducible reduced scheme of finite type over $\c$.

  The author expresses  her hearty thanks to Lawrence Ein  and Mircea  Musta\cedilla{t}\v{a} for the stimulating discussions.
  Before the publication of this paper, 
  Tommaso de Fernex and
  Roi Docampo introduced the same invariant and obtained the results same as our
Proposition \ref{mimic} and Corollary \ref{kata} independently (\cite{dfd}).
  
\section{Preliminaries on Mather discrepancy and arc spaces}

\noindent
Let $X$ be a $\q$-Gorenstein variety of index $r$ and $f: Y\to X$ a resolution of the singularities of $X$. Then the (usual) discrepancy divisor $K_{Y/X}$ is the unique 
  $\q$-divisor supported on the exceptional locus of $f$ such that $rK_{Y/X}$ is linearly equivalent with $rK_Y-f^*(rK_X)$.
 Note that the usual discrepancy is defined only for a $\q$-Gorenstein variety $X$, and 
 the following Mather discrepancy is defined for every variety, even for non-normal variety. 
\begin{defn} [\cite{DEI}]  Let $X$ be a variety of dimension $n$ and $f:Y\to X$ a resolution of the singularities factoring through the Nash blow up. 
Then, the image of the canonical homomorphism
$$f^*\wedge^n\Omega_X \to \wedge^n \Omega_Y$$
is an invertible sheaf of the form $J \wedge^n \Omega_Y$, 
where $J$ is the invertible ideal sheaf of $\o_Y$ that defines an effective divisor supported on the exceptional locus of $f$. 
This divisor is called the {\it Mather discrepancy divisor} and denoted by $\hKY$.
For every prime divisor $E$ on Y, we define
$$\widehat k_E:=\ord_E(\hKY). $$
More generally, if $v$ is a divisorial valuation over $X$, then we can assume without loss of generality that $v=q\val_E$ for a prime divisor $E$ on some $Y$ and a positive integer $q$, and define 
$$\widehat k_v:= q\cdot\widehat k_E.$$
For a $\q$-Gorenstein variety $X$, we define $k_E:=\ord_E(\KY)$ for a resolution $f:Y\to X$ and define also 
$k_v$ in the similar way.
\end{defn}

\begin{say}
\label{gorenst}
 Let $X$ be an $n$-dimensional $\q$-Gorenstein variety of index $r$.
 We write the image of the homomorphism
 $$(\wedge^n\Omega_X)^{\otimes r}\to \o_X(rK_X)=\omega_X^{[r]}$$
 by $I_r\otimes \o_X(rK_X)$, where $I_r$ is an ideal of $\o_X$.
 Let $f:Y\to X$ be a resolution factoring through the Nash blow up.
 Then, the relation of usual discrepancy and the Mather discrepancy is as follows:
 $$f^*(I_r)\otimes\o_Y(r\hKY)=\o_Y(r\KY).$$
 In particular $\hKY\geq\KY$.
 Let $\j_X$ be the Jacobian ideal of $X$ and let $J_r=(\overline{\j_X^r}:I_r)$, then $J_r\cdot I_r$ and $\j_X^r$ have the same integral closure (\cite[Corollary 9.4]{e-Mus2}).
 If $X$ is locally a complete intersection, then $I_1=\j_X$
\end{say}


\begin{defn}
  Let \( X \) be a scheme of finite type over \( \bC \)
and $K\supset \bC$ a field extension.
  For \( m\in \bN \), a  \( \bC \)-morphism \( \spec K[t]/(t^{m+1})\to X \) is called an \( m \)-jet
  of \( X \) and a  \( \bC \)-morphism \( \spec K[[t]]\to X \) is called an 
  {\it arc} of \( X \).
\end{defn}

\begin{say}
  We denote the space of $m$-jets of $X$ by $X_m$ and the space of arcs by 
  $X_\infty$. For  terminologies and the basic properties of these spaces,
  we refer the paper \cite{cr}.
\end{say}

\begin{defn}
   Let \( X \) be a variety over \( \bC\). 
   We say an arc \( \alpha:\sT  \to  X \) is {\it thin} 
   if \( \alpha  \) factors through a proper closed subset of \( X \).
  An arc which is not thin is called a {\it fat arc}.

  An irreducible closed subset \( C \) in \( X_{\infty} \)  is called a {\it 
  thin set} if  the generic point of \( C \) is thin.
  An irreducible closed subset  in \( X_{\infty} \) which is not thin  is called a {\it 
  fat set}.
  \end{defn}
  
  One typical example of a fat set is a maximal divisorial set which is introduced in \cite{i4}.

\begin{defn}
  For a divisorial valuation \( v \) over a variety \( X \), 
  define the maximal divisorial 
  set corresponding to \( v \) as follows:
  \[ C_{X}(v):=\overline{\{\alpha\in X_{\infty} \mid \ord_{\alpha}=v\}}, \]
  where \( \overline{\{ \ \ \}} \) is the Zariski closure in \( X_{\infty} \).
\end{defn}

\begin{prop}[\cite{i4}]
\label{irred}
  Let \( v=q\val_{E} \) be a divisorial valuation over a variety \( X \).
  Let \(f:Y\to X \) be a good resolution of the singularities of \( 
  X \) such that the prime divisor \( E \) appears on \( Y \).
  Here, a good resolution means a resolution whose exceptional locus is a 
  simple normal crossing divisor.
  Then, \[ C_{X}(v)=\overline{f_{\infty}(\cont ^q(E))}. \]
  In particular, \( C_{X}(v) \) is irreducible.  
\end{prop}

\begin{defn}[\cite{ELM}]
  For an ideal sheaf \( \a \) on a variety \( X \),
  we define 
   \[  \cont^e(\a)= \{\alpha \in X_{\infty} \mid \ord_{\alpha}(\a)= 
   e \} \] 
      and 
      \[  \cont^{\geq e}(\a)= \{\alpha \in X_{\infty} \mid \ord_{\alpha}(\a)
   \geq e \}. \] 
     These subset are called {\it contact loci} of an ideal \( \a \).
    The subset \( \cont^{\geq m}(\a) \) is closed and  
   \( \cont^m(\a) \) is locally closed.
  Both are cylinders.
 We can define in the obvious way also subsets $ \cont^e(\a)_m$ (if $e\leq m$) 
 and $\cont^{\geq e}(\a)_m$ (if $e\leq m+1$) in $X_m$.
\end{defn}

\begin{prop}[\cite{DEI}]
\label{component}
 Let \( X \) be an affine variety, and let \( \a_i\subset \o_X \) $(i=1,\ldots, r)$
be   non-zero ideals.
 Then, for \( e_1,\ldots, e_r\in \bN \),
 every fat irreducible component of  the  intersection \( \cont^{\geq
 {e_1}}(\a_1)\cap \cdots \cap \cont^{\geq
 {e_r}}(\a_r) \)  is a maximal divisorial set.
\end{prop}

Note that  \cite[Proposition 2.12]{DEI} is formulated for the case $r=1$.
But its proof  works also for $r>1$.

\begin{say}
  As $X$ is a variety over $\bC$, the arc space $X_\infty$ is irreducible (\cite{kln}).
  Therefore $\psi_m(X_\infty)$ is an irreducible constructible subset in $X_m$ of 
  dimension $(m+1)n$, where $n=\dim X$.
  Let $C\subset X_\infty$ be a cylinder $\psi_p^{-1}(A)$ contained in $\cont^e(\j_X)$.
  Then, codimension of $C$ is defined as follows:
  $$\codim (C,X_\infty):=(m+1)n-\dim \psi_m(C)$$
  for $m\geq \max\{p,e\}$.
  
  For an arbitrary cylinder $C$, the codimension is defined as follows:
  $$\codim (C,X_\infty):=\min \{\codim (C\cap \cont^e(\j_X))\mid e\in \bN\}.$$
  We sometimes write $\codim (C)$ for $\codim (C,X_\infty)$, when there is no
  possible confusion.
  Note that these are well defined by the following lemma. 
  (For details, see \cite[Section 5]{e-Mus2}.)
  \end{say}

\begin{lem}[\cite{DL2}, \cite{e-Mus2}]
\label{4.1}
Let $X$ be a variety of dimension $n$ and $e$ a nonnegative integer.
Let $\j_X$ be the Jacobian ideal of $X$.
Fix $m\geq e$ and let $\psi_{m+e, m}:X_{m+e}\to X_m$,
$\psi_m:X_\infty \to X_m$ be the truncation morphisms.
\begin{enumerate}
\item[(i)]
We have $\psi_m(\cont^e(\j_X))=\psi_{m+e, m}\left(\cont^e(\j_X)_{m+e}\right).$
\item[(ii)]
The truncation morphism $\psi_{m+1,m}:X_{m+1}\to X_m$ induces a piecewise trivial fibration
$$\psi_{m+1}\left(\cont^e(\j_X)\right)\to \psi_m\left(\cont^e(\j_X)\right)$$
with fiber $\bA^n$.
\end{enumerate}  
\end{lem}  

  \begin{lem}[\cite{e-Mus2}]
  \label{4.7}
  Let $X$ be a variety of dimension $n$. 
  If $m, p$ and $e$ are nonnegative integers such that $2p\geq m\geq p+e$,
  then the truncation morphism $\psi_{m,p}:X_m\to X_p$ induces a piecewise trivial
  fibration
  $$\cont^e(\j_X)_m\to \cont^e(\j_X)_p\cap \im(\psi_{m,p})$$
  with fiber $\bA^{(m-p)n+e}$.
    \end{lem}
    
    \begin{prop}[\cite{DEI}]
\label{cod}
  For a divisorial valuation $q\val_E$, we have the codimension of the corresponding maximal divisorial set
  $$q(\widehat k_E+1)=\codim (C_X(q\val_E)).$$
\end{prop}

\section{Invariants based on Mather discrepancy}
\noindent
First we start this section with the well known invariants.
\begin{defn} Let $(X,\a)$ be a pair consisting of $\q$-Gorenstein variety $X$ and a 
non-zero ideal $\a$ of $\o_X$.
The {\it log-canonical threshold} of $(X,\a)$ is defined as follows:
$$\lct(X,\a)=\sup\{c \mid k_E-c\ord_E(\a)+1 \geq 0, E\ \mbox{divisor\ over\ X}\}.$$

Let $W$ be a closed subset of $X$.
The {\it minimal log-discrepancy} of $(X,\a)$ along $W$ is defined as follows:

\noindent
If $\dim X\geq 2$,
$$\mld(W; X,\a)=\inf \{k_E-\ord_E(\a)+1 \mid E\ \mbox{divisor\ over}\ X \ \mbox{with\ center\ in\ } W\}.$$
When $\dim X=1$ we use the same definition of minimal log discrepancy, unless the infimum is negative,
in which case we make the convention that $\mld (W,X,\a)=-\infty$.
\end{defn}

\begin{rem}
\begin{enumerate}
\item[(i)]
The log-canonical threshold is also presented as 
$$\lct(X,\a)=\max\{c \mid k_{E_i}-c\ord_{E_i}(\a)+1\geq 0, \ E_i : \mbox{exceptional\ prime\ divisor\ on}\ Y\}$$
for a fixed log-resolution $f:Y\to X$ of $(X,\a).$
\item[(ii)]
If  $\mld(W;X,\a)<0$, then $\mld(W;X,\a)=-\infty$.
This is known when $\dim X\geq 2$, while it follows from the definition when $\dim X=1$. 
\end{enumerate}
\end{rem}

Now we will define the invariants modified from these invariants.
\begin{defn} Let $(X,\a)$ be a pair consisting of an arbitrary variety $X$ and a 
non-zero ideal $\a$ of $\o_X$.
The {\it Mather log-canonical threshold} of $(X,\a)$ is defined as follows:
$$\hlct(X,\a)=\sup\{c \mid \widehat k_E-c\ord_E(\a)+1 \geq 0, E\ \mbox{divisor\ over\ X}\}.$$

Let $W$ be a closed subset of $X$.
The {\it Mather minimal log-discrepancy} of $(X,\a)$ along $W$ is defined as follows:

\noindent
If $\dim X\geq 2$,
$$\hmld(W; X,\a)=\inf \{\ \widehat k_E-\ord_E(\a)+1 \mid E\ \mbox{divisor\ over}\ X
\ \mbox{with\ center\ in\ }W\}.$$
When $\dim X=1$ we use the same definition of Mather minimal log discrepancy, unless the infimum is negative,
in which case we make the convention that $\hmld (W,X,\a)=-\infty$.

\end{defn}

\begin{rem}
\label{attain}
\begin{enumerate}
\item[(i)]
  The Mather log-canonical threshold is represented as
$$\hlct(X,\a)=\max\{c \mid \widehat k_{E_i}-c\ord_{E_i}(\a)+1\geq 0, \ E_i : \mbox{exceptional\ prime\ divisor\ on}\ Y\}$$
for a fixed log-resolution $f:Y\to X$ of $(X,\a)$ factoring through the Nash blow up, because for a sequence 
$Y'\stackrel{g}\longrightarrow Y\stackrel{f}\longrightarrow X$ of such log resolutions 
of $(X,\a)$, we have $\widehat K_{Y'/X}=K_{Y'/Y}+g^*\hKY$ with $K_{Y'/Y}\geq 0$.
\item[(ii)]
If  $\hmld(W;X,\a)<0$, then $\hmld(W;X,\a)=-\infty$.
This is proved by using the previous formula of $\widehat K_{Y'/X}$ when $\dim X\geq 2$,
while it follows from the definition when $\dim X=1$.
\end{enumerate}
\end{rem}
\begin{prop} Let $X$ be an arbitrary variety and $\a$ is a non-zero ideal of $\o_X$.
Then,
$$\hlct(X, \a)=\min_{m\in \bN}\frac{\codim (\cont^{\geq m}(\a))}{m}.$$
\end{prop}

\begin{proof}
First we prove that
\begin{equation}
\label{leq}
\hlct(X, \a)\leq \inf _{m\in \bN}\frac{\codim (\cont^{\geq m}(\a))}{m}.
\end{equation}
For  $m\in \bN$, take an irreducible component $C_X(v)\subset\cont^{\geq m}(\a)$ which gives the codimension of $\cont^{\geq m}(\a)$.
Let $v=q\val_E$ for a prime divisor $E$ over $X$ and let $r_E:=\ord_E(\a)$.
Then, $C_X(v)\subset \cont^{\geq m}(\a)$ implies $qr_E=v(\a)\geq m.$
Therefore,
$$\codim(\cont^{\geq m}(\a))=\codim(C_X(v))=q(\widehat k_E+1)
\geq m\frac{\widehat k_E+1}{r_E}.$$
Here, by the definition of $\hlct$, the last term is  $\geq m\cdot\hlct(X,\a)$, which gives the inequality (\ref{leq}).
For the theorem, it is sufficient to prove that there is $m$ such that
\begin{equation}
\label{geq}
 \frac{\codim (\cont^{\geq m}(\a))}{m}\leq \hlct(X, \a).
\end{equation}
By Remark \ref{attain} (i), there is a prime divisor $E$ over $X$ such that 
$\hlct(X,\a)={(\widehat k_E+1)}/{\ord_E(\a)}$.
Let $m={\ord_E(\a)}$, then we have the inclusion $C_X(\val_E)\subset \cont^{\geq m}(\a)$  and this implies
the inequality
$$\codim (\cont^{\geq m}(\a))\leq \codim(C_X(\val_E))=\widehat k_E+1=m\cdot \hlct (X,\a),$$ which gives the required inequality (\ref{geq}).
\end{proof}

As a corollary, we obtain  the formula of $\lct$ for non-singular case. 

 \begin{cor}[\cite{ELM}]
 \label{smoothlct}
 Let  $(X, \a)$ be a pair consisting of a non-singular variety $X$  and an ideal 
 $\a\subset \o_X$. Let $Z$ be the subscheme defined by $\a$.
 Then the log-canonical threshold is obtained as follows:
 $$\lct(X,\a)=\min_{m\in \n}\frac{\codim(Z_{m-1},X_{m-1})}{m}.$$
\end{cor}

This follows immediately from the theorem, since the equality \\
$\codim (\cont^{\geq m}(\a))=\codim(Z_{m-1}, X_{m-1})$ holds for non-singular $X$.

The next is the formula for the Mather minimal log-discrepancy in terms of the arc space.

\begin{prop}
\label{mldformula}
 Let  $(X, \a)$ be a pair consisting of an arbitrary variety $X$  and a non-zero ideal 
 $\a\subset \o_X$.
Let $W$ be a proper closed subset of $X$ and $I_W$ be the (reduced) ideal of $W$.
Then,
\begin{equation}
\label{hmld=}
\hmld(W;X,\a)=\inf_{m\in\n}\{\codim(\cont^m(\a)\cap\cont^{\geq 1}(I_W))-m\}.
\end{equation}
We also have 
\begin{equation}
\label{hmldgeq}
\hmld(W;X,\a)=\inf_{m\in\n}\{\codim(\cont^{\geq m}(\a)\cap\cont^{\geq 1}(I_W))-m\}.
\end{equation}

\end{prop}

\begin{proof}
  For the proof of $\geq$ at (\ref{hmld=}), let $E$ be any prime divisor over 
  $X$ with the center in $W$. 
  Let $m=\ord_E(\a)$ and $v=\val_E$.
  Then, there is a non-empty open subset $C$ of $C_X(v)$ such that $C\subset \cont^m(\a)\cap\cont^{\geq 1}(I_W)$.
  Hence, 
    $$\widehat k_E-\ord_E(\a)+1=\codim(C_X(v))-m$$
   $$\geq 
    \codim(\cont^m(\a)\cap\cont^{\geq 1}(I_W))-m,$$
    which yields the required inequality unless $\dim X=1$ and $\hmld(X,\a)=-\infty$.
    
    When $\dim X=1$ and $\hmld(X,\a)=-\infty$, there is a prime divisor $E$ over $X$ with the center in $W$
    such that $\widehat k_E-\ord_E(\a)+1<0$.
    Let $m=\ord_E(\a)$, then 
    $\codim C_X(\val_E)-m<0.$
    Here, for every $q\in \bN$, by Proposition \ref{cod},
    $$\codim C_X(q\val_E)-qm=q(\codim C_X(\val_E)-m)<0.$$
    As a non-empty open subset of $C_X(q\val_E)$ is in $\cont^{ qm}\cap\cont^{\geq 1}(I_W)$,
    we have 
    $$\codim(\cont^{ qm}\cap\cont^{\geq 1}(I_W))-qm\leq \codim (C_X(q\val_E))-qm$$ 
    $$=q(\codim (C_X(\val_E))-m)<0.$$ 
    Here, if $q\to \infty$, then we have $$\codim(\cont^{qm}\cap\cont^{\geq 1}(I_W),)-qm
    \to -\infty,$$ which implies 
   the right hand side of (\ref{hmld=}) in the theorem is $-\infty$.

    For the proof of $\leq$ at (\ref{hmld=}),
    we may assume that $\widehat k_E-\ord_E(\a)+1\geq 0$
    for every prime divisor $E$ over $X$ with the center in $W$.
    Indeed if there is a prime divisor $E$ with the  center in $W$ and $\widehat k_E-\ord_E(\a)+1<0$, 
    then $\hmld(W;X,\a)=-\infty$ by Remark \ref{attain}, (ii) and therefore the required 
    inequality is trivial.
    
    For $m\in \bN$, let $C\subset 
    \cont^m(\a)\cap\cont^{\geq 1}(I_W)$ be an irreducible component
    which gives the codimension of $\cont^m(\a)\cap\cont^{\geq 1}(I_W)$.
    Then, the closure $\overline C$ is $C_X(v)$ for some  divisorial valuation $v$,
    since a fat irreducible component of $\cont^{\geq m}(\a)\cap\cont^{\geq 1}(I_W)$
    is a maximal divisorial set (Proposition \ref{component}) and 
    $\cont^m(\a)\cap\cont^{\geq 1}(I_W)$ is an open subset of 
    $\cont^{\geq m}(\a)\cap\cont^{\geq 1}(I_W)$.
    Let $v=q\val_E$ and $m=v(\a)$, then $E$ is a prime divisor over $X$ with the center in $W$, $m=q\ord_E(\a)$ and
    $$\codim(\cont^m(\a)\cap\cont^{\geq 1}(I_W))-m
    =\codim(C_X(v))-m$$
    $$=	q(\widehat k_E+1)-q\ord_E(\a)\geq \widehat k_E+1-\ord_E(\a),$$
    which yields the required inequality.
    
    For the proof of (\ref{hmldgeq}) of the theorem, let
    $$a_m=\codim(\cont^m(\a)\cap\cont^{\geq 1}(I_W))-m,$$ 
$$b_m=\codim(\cont^{\geq m}(\a)\cap\cont^{\geq 1}(I_W))-m.$$
  As $\cont^m(\a)\subset \cont^{\geq m}(\a)$, we have $a_m\geq b_m$.
  Therefore, it follows\\ $\inf_m \{a_m\}\geq \inf _m\{b_m\}$.
  
  Next we prove the converse inequality.
  For every $m\in \bN$, let $C_X(v)$ be the irreducible component of 
  $\cont^{\geq m}(\a)\cap\cont^{\geq 1}(I_W)$ that gives the codimension.
  Then, for $m':=v(\a)\geq m$ we have 
  $$\codim(\cont^{\geq m}(\a)\cap\cont^{\geq 1}(I_W))=\codim(\cont^{m'}(\a)\cap\cont^{\geq 1}(I_W)).$$
  Hence, $b_m\geq a_{m'}$, which yields $\inf_m \{b_m\}\geq \inf _m\{a_m\}$.
\end{proof}

\begin{rem}
\label{extended}
 Our formula can be easily extended for the combination of ideals 
 $\a_1,\a_2,\cdots,\a_r$ instead of one ideal $\a$. I.e., we have
 $$\hmld(W;X,\a_1^{e_1}\a_2^{e_2}\cdots\a_r^{e_r})=$$
 $$\inf_{m_i\in\n}\{\codim(\cont^{ m_1}(\a_1)\cap\cdots\cap\cont^{ m_r}(\a_r)\cap\cont^{\geq 1}(I_W))-\sum_im_ie_i\},$$
 where $e_i$'s are positive real numbers.
 Here, any of $\cont^{m_i}(\a_i)$'s can be replaced by $\cont^{\geq m_i}(\a_i)$.
 For  simplicity of the notation and the proofs, we keep formulating the forthcoming formulas 
 for one ideal only.
 But note that the formulas in this section are also valid under this combination form.

\end{rem}

The following  lemma is a modified version of  \cite[Theorem 7.9]{e-Mus2} for our discussion.

\begin{lem} 
\label{28}
Let $X$ be a variety of dimension $n$, then
$$\hmld(W;X, \a\j_X)=$$
$$\inf \left\{(s+1)n-m-\dim \left(\cont^{\geq m}(\a)_s\cap\cont^e(\j_X)_s\cap \cont ^{\geq 1}(I_W)_s \right)\right\},$$
where the infimum is over those $m,e\in \n$ and $s\in \n$ such that $s\geq 2e, e+m$.
\end{lem}

\begin{proof}
  Let $S:=\cont^{\geq m}(\a)_s\cap\cont^e(\j_X)_s\cap\cont^{\geq 1}(I_W)_s\subset X_s$.
  Then, by Lemma \ref{4.7}, the truncation morphism $ \psi_{s,s-e}$ induces
  a piecewise trivial fibration
  $S\to \psi_{s,s-e}(S)$
  with fiber $\bA^{e(n+1)}$.
 By Lemma \ref{4.1} (i), it follows 
 $$\psi_{s,s-e}(S)=\psi_{s-e}(\cont^{\geq m}
 (\a)\cap\cont^e(\j_X)\cap\cont^{\geq 1}(I_W)).$$
 On the other hand, by Lemma \ref{4.1} (ii), 
 $$\codim (\cont^{\geq m}(\a)\cap\cont^e(\j_X)\cap\cont^{\geq 1}(I_W))$$
 $$=(s-e+1)n-\dim \left(\psi_{s-e}(\cont^{\geq m}(\a)\cap\cont^e(\j_X)\cap\cont^{\geq 1}(I_W)\right),$$
 as $s-e\geq m, e, 1$.
 The last term is equal to $(s-e+1)n-(\dim S -e(n+1))$, by the previous discussion.
 Now we obtain 
 $$
 \begin{array}{ll}
 \hmld(W;X,\a\j_X)&=\displaystyle{\inf_{m,e}}\{\codim (\cont^{\geq m}(\a)\cap\cont^e(\j_X)\cap\cont^{\geq 1}(I_W))-m-e\}\\
&=\displaystyle{\inf}\{(s+1)n-m-\dim S\},
\end{array}$$
where the infimum is over those $m,e\in \n$ and $s\in \n$ such that $s\geq 2e, e+m$.
\end{proof}

The following is a generalization of \cite[Theorem 8.1]{e-Mus2} and \cite{kaw1}, 
whose formulation is in Corollary \ref{otehon}.
But our proof is just an imitation of the proof of \cite[Theorem 8.1]{e-Mus2} and  is 
even easier.

\begin{prop}[Inversion of Adjunction]
\label{mimic}
Let $X$ be an arbitrary varity, $A$ a non-singular variety containing $X$ as a closed
subvariety of codimension $c$ and $W$ a proper closed subset of $X$.
 Let $\widetilde\a\subset \o_A$ be an ideal such that its image ${\a}:=\widetilde\a\o_X\subset \o_X$ is non-zero.
 Denote the ideal of $X$ in $A$ by $I_X$.
Then,
$$\hmld(W; X,{\a}\j_X)=\hmld(W;A,\widetilde\a I_X^c).$$
\end{prop}

\begin{proof}
  Our proof is based on the proof of  \cite[Theorem 8.1]{e-Mus2} which is a
formula of usual mld for a  
$\q$-Gorenstein variety $X$. 
  For  the reader's convenience, here we write down the proof with a care on 
  the difference between  their statement and ours.  
  
  For the proof of $\geq$ in the equality of the proposition, let 
  $\sigma=\hmld(W;A,\widetilde\a I_X^c)$ and let $n=\dim X$.
  We assume contrary, $\hmld(W; X,{\a}\j_X)<\sigma$, and  will induce a 
  contradiction.
  By the assumption and Lemma \ref{28}, there exist $e, m\in \bN$, $s\geq 2e, m+e$ 
  and an irreducible component $C$ of $\cont^{\geq m}(\a)_s\cap\cont^e(\j_X)_s\cap\cont^{\geq 1}(I_W)_s$ such that 
\begin{equation}
\label{small}
(s+1)n-m-\dim C<\sigma,
\end{equation}
where $I_W$ is the defining ideal of $W$ in $X$.
As $C\subset \cont^{\geq m}(\a)_s\cap\cont^{\geq 1}(I_W)_s=
\cont^{\geq m}(\widetilde\a)_s\cap\cont^{\geq 1}(\widetilde I_W)_s \cap X_s,$
where $\widetilde I_W$ is the defining ideal of $W$ in $\o_A$,
we have 
$$(\psi_s^A)^{-1}(C)\subset  
\cont^{\geq m}(\widetilde\a)\cap\cont^{\geq 1}(\widetilde I_W)\cap\cont^{\geq s+1}(I_X)=:S,$$
where $\psi^A_s: A_\infty \to A_s$ is the truncation morphism. 
Here we note that $\cont^{\geq s+1}(I_X)=\cont^{\geq c(s+1)}(I_X^c)$.
Now we obtain 
$$(n+c)(s+1)-\dim C=\codim (C, A_s)=\codim ((\psi_s^A)^{-1}(C), A_\infty)$$
$$\geq \codim (S, A_\infty) \geq \sigma +m+c(s+1),$$
which is a contradiction to (\ref{small}).
 
 For the proof of $\leq$ in the equality of the proposition, let $\tau=
 \hmld(W; X,{\a}\j_X)$.
 We assume contrary, $\tau> \hmld(W;A,\widetilde\a I_X^c)$, and will induce a contradiction.
 By the assumption and Proposition \ref{mldformula}, there exists an irreducible 
 component $C$ of $\cont^{\geq m}(\widetilde \a)\cap\cont ^{\geq d}(I_X)\cap
 \cont^{\geq 1}(\tilde I_W)$ such that 
 \begin{equation}
 \label{big}
 \codim (C, A_\infty)<cd+m+\tau.
 \end{equation}
 Just as in the proof of \cite[Theorem 8.1]{e-Mus2}, there is an open subcylinder 
 $C_0\subset C$ such that $\codim (C_0\cap X_\infty, X_\infty)\leq \codim (C_0, A_\infty)+
 e-cd,$ where $e=\min_{\gamma\in C\cap X_\infty}\{\ord_\gamma(\j_X)\}$.
 Therefore, we have 
 $$\codim (C_0\cap X_\infty, X_\infty)< m+\tau+e.$$
 On the other hand, as $C_0\cap X_\infty\subset \cont^{\geq m}(\a)\cap \cont^{\geq e}(\j_X)\cap\cont^{\geq 1}(I_W)$, we obtain 
  $$\tau \leq \codim( \cont^{\geq m}(\a)\cap \cont^{\geq e}(\j_X)\cap\cont^{\geq 1}(I_W))-m-e$$ 
 $$\leq \codim (C_0\cap X_\infty, X_\infty)-m-e<\tau,$$
 which is a contradiction.
\end{proof}

As 
$\mld(W; X, \a J_r^{1/r})=\hmld(W;X,\a\j_X)$ for $\q$-Gorenstein $X$
we obtain the following:
\begin{cor}[\cite{EMY}, \cite{e-Mus2}, \cite{kaw1}]
\label{otehon}
Let  $X$ be a normal closed subvariety in a non-singular  variety $A$ of codimension $c$ and
 let $W$ be a proper closed subset of $X$. 
Assume that $X$ is $\q$-Gorenstein variety of index $r$.
Let $\widetilde\a\subset \o_A$  be an ideal such that its image  ${\a}:=\widetilde\a\o_X\subset \o_X$ is non-zero. Then,
$$\mld(W;X,\a J_r^{1/r})=\mld(W; A, \widetilde\a I_X^c),$$
where $I_X$ is the defining ideal of $X$ in $A$ and $J_r$ is as in \ref{gorenst}.
\end{cor}

\begin{cor}[Adjunction formula]
\label{nonequal}
  Let $X$ be a closed subvariety of a variety $X'$ of codimension $c$ and let $W$ be a proper closed subset of $X$. Let $\a'\subset \o_{X'}$ be an ideal such that its image 
  $\a:=\a'\o_X \subset \o_X$ is non-zero. Let $I_{X/{X'}}$ be the defining ideal of $X$ in 
  $X'$.  Then,
  $$\hmld(W; X,{\a}\j_X)\geq \hmld(W;X',\a'\j_{X'} I_{X/{X'}}^c).$$
\end{cor}

\begin{proof}
  We may assume that $X'$ is affine. Let $A$ be a non-singular variety containing 
  $X'$ and $I_{X/A}$ and $I_{X'/A}$ be the ideals of $X$ and $X'$ in $A$, respectively. Let $\widetilde\a\subset \o_A$ is an ideal whose image in $\o_{X'}$ is 
  $\a'$.
  Let $c'$ be the codimension of $X'$ in $A$.
  By Proposition \ref{mimic}, we have
  $$\hmld(W; X,{\a}\j_X)=\hmld(W;A,\widetilde\a I_{X/A}^{c+c'}), \mbox{and}$$
  $$\hmld(W; X',{\a'}\j_{X'}I_{X/{X'}}^c)=\hmld(W;A,\widetilde\a I_{X/A}^cI_{X'/A}^{c'}).$$
  As we have an inclusion 
  $$\widetilde\a I_{X/A}^{c+c'}\supset \widetilde\a I_{X/A}^cI_{X'/A}^{c'}$$
  of ideals,  the inequality of our corollary follows.
\end{proof}

\begin{exmp}
[L. Ein, M. Musta\cedilla{t}\v{a}]
The inequality in Corollary \ref{nonequal} is not the equality in general.
Indeed the following is a counter example.
Let $A=\bA^4_\bC$, $X'\subset A$ be a hypersurface defined by $xy-zw=0$ and 
$X\subset A$ be defined by $x=z=0$.
Let $0$ be the origin of $A$.
Then,  $$\hmld(0; X,\j_X)=2>1\geq \hmld(0; X',\j_{X'}I_{X/{X'}}).$$
Actually, as $X$ is a non-singular surface, $\hmld(0; X,\j_X)=\mld(0;X,\o_X)=2$.
On the other hand, by Proposition \ref{mimic}, it follows
$$\hmld(0; X',\j_{X'}I_{X/{X'}})=\mld(0; A, I_{X/A}I_{X'/A}).$$
Let $E$ be the exceptional divisor of the blow-up $\widetilde A\to A$ at the origin $0$.
Then, $\ord_EI_{X/A}I_{X'/A}=3$, since the ideal is a homogeneous ideal generated by polynomials of degree 3.
Hence,
$$\ord_EK_{{\tilde A}/A}-\ord_EI_{X/A}I_{X'/A}+1=1,$$
which yields $\mld(0; A, I_{X/A}I_{X'/A})\leq 1$.

\end{exmp}
\begin{cor} Let $(X,\a)$ be a pair consisting of an arbitrary variety $X$ and a non-zero ideal 
$\a\subset \o_X$, then the function $x\mapsto
\hmld(x;X,\a\j_X) $, ($x\in X$ closed point) is lower semicontinuous.
\end{cor}

\begin{proof} The following proof is the same as in \cite{e-Mus}.
  Since the statement is local, it is sufficient to prove the corollary for an affine variety 
  $X$. 
  So, we may assume that $X$ is a closed subvariety of an affine space $A=\bA^N$.
  By Proposition \ref{mimic}, for $x\in X$,
  $$\hmld(x; X,{\a}\j_X)=\hmld(x;A,\widetilde\a I_X^c).$$
  Here, the right hand side is $\mld(x;A,\widetilde\a I_X^c)$ since $A$ is non-singular and it is known that the function $x\mapsto
\mld(x;A,\widetilde\a I_X^c) $ is lower semicontinuous for non-singular $A$ (\cite{florin2}, \cite{EMY}).
Now we have only to restrict the function on $X$.
\end{proof}
  In general, the map $x\mapsto
\hmld(x;X,\a) $, ($x\in X$ closed point) is not necessarily lower semi-continuous (see 
Example \ref{non-lower}).

\begin{cor}
\label{kata}
 Let $X$ be a variety of dimension $n$. Then, for every closed point 
$x\in X$, the following inequality holds:
$$\hmld(x; X, \j_X)\leq n,$$
where the equality holds if and only if $(X,x)$ is non-singular.
\end{cor}

\begin{proof}
The first statement is an immediate consequence of the lower semicontinuity.
But for the second statement we give a different and direct proof of the first 
statement.
We may assume that $X$ is affine, therefore it is a closed subvariety of an affine 
space $A=\bA^N$. 
Let $I_X$ and $M_x$ be the defining ideals of $X$ and $x$ in $A$ and $c=N-n$.
  By Proposition \ref{mimic},
  $$\hmld(x; X, \j_X)=\hmld(x;A,I_X^c).$$
  By Proposition \ref{mldformula} and Remark \ref{extended}
 $$\hmld(x;A,I_X^c)=\inf_m\{\codim (\cont^{\geq m}(I_X)\cap \cont^{\geq 1}(M_x), A_\infty)-cm\}$$
$$=\inf_m\{ \codim(X_{m-1}\cap \pi_{m-1}^{-1}(x), A_{m-1})-cm\},$$
 where $\pi_{m-1}:A_{m-1}\to A$ is the canonical projection.
 Since the restiriction   $\overline{(X_{reg})_{m-1}}\to X$ of the projection has the relative dimension $(m-1)n$, it follows
 \begin{equation}
 \label{shokurov}
(m-1)n\leq \dim \overline{(X_{reg})_{m-1}}\cap \pi_{m-1}^{-1}(x)\leq
 \dim X_{m-1}\cap \pi_{m-1}^{-1}(x).
 \end{equation}
 Therefore, for every $m\in \bN$,
 $$ \hmld(x;A,I_X^c)\leq \dim A_{m-1}-(m-1)n-cm=n.
$$
This shows the first statement.

If $(X,x)$ is non-singular, then $\hmld(x;X,\j_X)=\mld(x;X,\o_X)$ and it is well known that $\mld(x;X,\o_X)=n$.
Conversely, if $\hmld(x;X,\j_X)=n$, then the equalities in (\ref{shokurov}) hold for every $m\in \bN$.
Hence, in particular for $m=2$, 
we have
$\dim X_1\cap\pi_1^{-1}(x)=n$.
This yields $\dim T_{X,x}=n$, which means that $(X,x)$ is non-singular.
\end{proof}

In \cite{katata} Shokurov posed the following conjecture:
\begin{conj}
  Let $X$ be a $\q$-Gorenstein variety of dimension $n$. Then, for every closed point 
$x\in X$, the following inequality holds:
$$\mld(x; X, \o_X)\leq n,$$
where the equality holds if and only if $(X,x)$ is non-singular.
\end{conj}

Our Corollary \ref{kata} is the answer to a modified version of this conjecture. 
In particular, if $(X,x)$ is a complete intersection, then the affirmative answer to 
this conjecture follows from our corollary, because 
$\mld(x; X, \o_X)=\hmld(x; X, \j_X)$.
This is already observed by Florin Ambro (private communication to the author) who published its special case  in \cite{florin}.

\section{A characterization of the singularities with the minimal Mather discrepancy}

\noindent
In this section we think of only the pair $(X,\o_X)$ for  a variety $X$ of dimension $n\geq 2$. 
We write $\hmld (x; X)$ for $\hmld (x; X,\o_X)$.
We denote the canonical projections of jet schemes by $\psi_m:X_\infty\to X_m$,
$\psi_{m',m}:X_{m'}\to X_m \ (m'>m)$, $\pi_m:X_m\to X$
and $\pi:X_\infty\to X$.
When we denote the projections for a variety, say $Y$, different from $X$, we put the symbol of 
the variety on the shoulder like $\psi^Y_m$ etc..

\begin{defn}
  Let $x$ be a closed point of a variety $X$ of dimension $n$.
  For every $m\in \bN$, define $\lambda_m$ and $\lambda^0_m$ as follows:
  $$\dim \psi_m(\pi^{-1}(x))=mn-\lambda_m,$$
  $$\dim \psi_m(\pi^{-1}(x)\setminus (\sing X)_\infty)=mn-\lambda^0_m,$$
  where $\sing X$ is the singular locus of $X$.
\end{defn}

The invariant $\lambda_m$ is introduced by Lawrence Ein.
Note that $\lambda_m=\lambda^0_m$ $(m\in \bN)$ if $(X,x)$ is an isolated singularity by \cite[Lemma 2.12]{i-k}.
(This is also proved  by the fact that $(\sing X)_\infty$ is consisting of the trivial arc on $x$.)
We have some basic properties of these invariants as  follows:

\begin{lem}
\label{lambda}
\begin{enumerate}
\item[(i)]
 For every $m\in \bN$, it follows $\lambda^0_m\geq \lambda_m\geq 0$.
\item[(ii)]
For $m\gg 0$, it holds that $\lambda_m=\lambda^0_m=\lambda$ (constant).
 \item[(iii)]
 $\lambda=\hmld(x;X)-n$. 
\end{enumerate}
\end{lem}

\begin{proof}
By the definition, it is clear that $\lambda_m\leq \lambda^0_m$.
By Denef-Loeser \cite[Lemma 4.3]{DL2}, the canonical projection $\psi_{m,l}$ 
induces a map $\psi_{m}(X_\infty)\to \psi_l(X_\infty)$ whose fibers have dimension $\leq (m-l)n$, where $m>l$ are non-negative integers,
it follows that $\psi_m(\pi^{-1}(x))=\psi_m(X_\infty)\cap \pi_m^{-1}(x)$ is of dimension $\leq mn$ by thinking of $l=0$. This yields (i).

For (ii), remember that $\pi^{-1}(x)$ is a contact locus and therefore has a finite number of irreducible components $C_1,\ldots, C_r, Z_1,\ldots,Z_s$, 
where $C_i$'s are fat, while $Z_j$'s are thin and contained in $(\sing X)_\infty$
(\cite{DEI}).
As $Z_j$ is thin, we have $$\dim \psi_m(Z_j)\leq m(n-1)$$ for every $m\in \bN $ and $j$ by the discussion on (i).
On the other hand, by Proposition \ref{cod}, for $m\gg 0$,
$$\dim \psi_m(C_i)=(m+1)n-(\widehat k_{E_i}+1),$$
where $C_i=C_X(\val_{E_i})$. 
(Note that the valuation is reduced, since $C_i$ is a component of $\pi^{-1}(x)$.)
Let $\lambda=\min_i\{\widehat k_{E_i}+1\}-n$, then for $m\gg 0$ we have 
$mn-\lambda_m=\dim \psi_m(\pi^{-1}(x))=mn-\lambda$.
For $\lambda^0_m$, we have the required form from 
$$\dim \psi_m(\pi^{-1}(x)\setminus (\sing X)_\infty)= \dim \psi_m(C_1\cup \cdots \cup
C_r)$$
for $m\gg 0$.

For (iii), we have only to note that $$\pi^{-1}(x)=\bigcup_{v}C_X(v),$$
up to a thin set, where $v $ are all divisorial valuation centered at $\{x\}$.
Thus, $\lambda=\hmld(x;X)-n$.
\end{proof}

The following is a trivial consequence from the lemma.
\begin{cor} For every closed point $x$ of a variety $X$ of dimension $n$,
$$\hmld(x;X)\geq n.$$
\end{cor}

For the characterization of singularities with $\hmld(x;X)= n$, we prepare a few lemmas.

\begin{lem}
\label{reg}
  Let $b:Y\to X$ be the blow-up at the point $x$ and $g:\tilde Y\to Y$ be a resolution of the singularities of $Y$.
  Let $h:\tilde Y\to X$ be the composite $b\circ g$.
   Let $\pi_1^{\tilde Y}:\tilde Y_1\to \tilde Y $ be the canonical projection from 1-jet scheme 
   $\tilde Y_1$ of $\tilde Y$ to $\tilde Y$.
  Denote the scheme theoretic fibers $b^{-1}(x)$ and $h^{-1}(x)$ by $E$ and $E'$,
  respectively.
  Then,
  $$\overline{\psi_1\left(\pi^{-1}(x)\setminus (\sing X)_\infty\right)}=
  \overline{h_1\left((\pi_1^{\tilde Y})^{-1}(E')\right)},$$
  and if this subset contains a non-trivial jet, then we also have the following equality:
  $$ \overline{h_1\left((\pi_1^{\tilde Y})^{-1}(E')\right)}=
 \overline{h_1\left((\pi_1^{\tilde Y})^{-1}(E'_{reg})\right)},$$
  where $E'_{reg}$ is the non-singular locus of the scheme $E'$.
\end{lem}

\begin{proof}
  As $h$ is a proper morphism, by the valuation criteria of properness, 
  every arc in $\pi^{-1}(x)\setminus (\sing X)_\infty$ can be lifted to an arc
  in $(\pi^{\tilde Y})^{-1}(E')\setminus (h^{-1}(\sing X))_\infty$.
  (When $x$ is a non-singular point, just use the surjectivity of 
  $h_\infty:\tilde Y_\infty \to X_\infty$ on a neighborhood  of $x$.)
  Therefore,
  $$\psi_1(\pi^{-1}(x)\setminus (\sing X)_\infty)=h_1\psi_1^{\tilde Y}\left((\pi^{\tilde Y})^{-1}(E')\setminus (h^{-1}(\sing X))_\infty\right).$$
  Noting that 
$(\pi^{\tilde Y})^{-1}(E')\setminus (h^{-1}(\sing X))_\infty$ is dense in $(\pi^{\tilde Y})^{-1}(E')$, we have the first equality of the lemma.
As for the second equality, the inclusion $\supset$ is trivial. 
To show the inclusion $\subset$, 
take any non-trivial jet $h_1(\alpha)\in 
h_1(\pi_1^{\tilde Y})^{-1}(E')$ for $\alpha\in ( \pi_1^{\tilde Y})^{-1}(E')$.
Let $P=\pi_1^{\tilde Y}(\alpha)\in E'$.
Then $\alpha$ corresponds to a ring homomorphism 
$$\alpha^*:\o_{\tilde Y,P}\to \c[t]/(t^2)$$
which sends the maximal ideal ${\frak m}_{\tilde Y, P}$ to $(t)$.
If the homomorphism $h^*:\o_{X,x}\to \o_{\tilde Y, P}$ corresponding to $h$ 
sends ${\frak m}_{X,x}$ to ${\frak m}^2_{\tilde Y, P}$, then $\alpha^*h^*({\frak m}_{X,x})=0$ which means that $h_1(\alpha)$ is the trivial 1-jet, a contradiction to 
the definition of $h_1(\alpha)$.
Therefore the principal ideal ${\frak m}_{X.x}\o_{\tilde Y,P}\not\subset
{\frak m}^2_{\tilde Y, P}$, which yields that $E'$ is non-singular at $P$.
Now we have $$\alpha\in ( \pi_1^{\tilde Y})^{-1}(E'_{reg}).$$
On the other hand,  the trivial jet is contained in the  subsets $$ \overline{{h_1\left((\pi_1^{\tilde Y})^{-1}(E')\right)}\setminus 
\{\mbox{the\ trivial\ jet}\}},$$ by the previous argument, the trivial jet is contained in
$ \overline{h_1\left((\pi_1^{\tilde Y})^{-1}(E'_{reg})\right)}$, which completes the proof of the second equality
\end{proof}

\begin{say}
 Let $T_{X,x}$ be the tangent space of $X$ at $x$ and $C_{X,x}$  the tangent cone of $X$ at $x$.
 Let $b:Y\to X$ be the blow-up of $X$ at $x$ and $E=b^{-1}$(x) the scheme theoretic fiber.
 Then, it is well known that  there is an isomorphism 
 $$E\simeq \left(C_{X,x}\setminus \{0\}\right)/\c^* \subset \left(T_{X,x}\setminus \{0\}\right)/\c^* .$$
Denote the projection $T_{X,x}\setminus \{0\}\to T_{X,x}\setminus \{0\}/\c^*$ by $\rho$.

\end{say}

\begin{lem}
\label{corres}
 Let $\pi_1:X_1\to X$ be the canonical projection from the 1-jet scheme.
Then, we have the following:
\begin{enumerate}
\item[(i)]
For a closed point $x\in X$, there is an isomorphism 
$$\phi:\pi_1^{-1}(x)\isom T_{X,x},$$ 
which sends the trivial jet to the origin $0$ 
and  sends the subset $\psi_1\left(\pi^{-1}(x)\right)$ into 
$C_{X,x}$.
\item[(ii)] For a point $P\in E$, take a 1-jet $\alpha\in (\pi_1^Y)^{-1}(P)$ of $Y$.
  If $b_1(\alpha)$ is a non-trivial jet, then, the point 
  $\rho\circ\phi(b_1(\alpha))\in \left(C_{X,x}\setminus \{0\}\right)/\c^*$ corresponds to
  $P$ under the isomorphism $\left(C_{X,x}\setminus \{0\}\right)/\c^*\isom E$.
%
  \end{enumerate}
\end{lem}

\begin{proof}
(i) The existence of the isomorphism $\phi$ is well known by the definitions of $\pi_1^{-1}(x)$ and Zariski tangent space.
Here, we give a concrete description of $\phi$ for the rest of the statements (i).

Let $(X,x)\subset (\bA^N,0)$ and let $I$ be the defining ideal of $X$ in $\bA^N$.
A jet $\alpha\in \pi_1^{-1}(x)$ corresponds to a ring homomorphism 
$$\alpha^*:\c[x_1,\ldots,x_N]\to \c[t]/(t^2)\ ;\ x_i\mapsto a_i^{(1)}t,$$
with $\alpha^*(I)=0$.
This homomorphism corresponds exactly to the point $$(a_1^{(1)},\ldots, a_N^{(1)})\in Z(I_1)\subset \bA^N,$$
where $I_1$ is the ideal of the linear terms of $I$.
Here, $Z(I_1)$ is the tangent space.
This correspondence $\alpha \mapsto (a_1^{(1)},\ldots, a_N^{(1)})$ is $\phi$.
It is clear that the trivial jet corresponds to the origin.
If a $1$-jet $\alpha$ can be lifted to an arc $\tilde \alpha$, then there is a commutative diagram of ring homomorphisms:
$$
\xymatrix
{& \c[[t]]\ar[d]\\
{ \c[x_1,\ldots,x_N]}\ar[ru]^{\tilde\alpha^*}\ar[r]^{\alpha^*}&\c[t]/(t^2)\\}
$$

with $\tilde\alpha^*(I)=0$, which yields $f_*(a_1^{(1)},\ldots,a_N^{(1)})=0$ for the initial term $f_*$ of every element $f\in I$. 
Therefore, this jet corresponds to a point 
$$(a_1^{(1)},\ldots, a_N^{(1)})\in Z(I_*)\subset \bA^N,$$
where $I_*$ is the ideal of the initial terms of $I$.
Here, $Z(I_*)$ is the tangent cone.
This completes the proof of (i).

For the proof of (ii), we may assume that $(X,x)=(\bA^N,0)$.
The blow-up $Y$ is covered by 
$U_{(i)}=\spec \c[x_i,\frac{x_1}{x_i},\ldots, \frac{x_N}{x_i}] $ $(i=1,\ldots, N)$.
We may assume that $P\in U_{(1)}$.
Then, $\alpha$ corresponds to the ring homomorphism
$$\alpha^*:\c\left[x_1,\frac{x_2}{x_1},\ldots, \frac{x_N}{x_1}\right] \to \c[t]/(t^2)$$
$$x_1\mapsto a_1^{(1)}t,\ \  \frac{x_i}{x_1}\mapsto  a_i^{(0)}+ a_i^{(1)}t, (i\neq 1).$$
Here, we note that the homogeneous coordinates of $P\in E\simeq \bP^{N-1}$
is $$(1:a_2^{(0)}:\ldots:a_N^{(0)}).$$
As $b_1(\alpha)$ corresponds to the ring homomorphism 
$$b_1(\alpha)^*=\alpha^*\circ b^*:\c[x_1,\ldots, x_N]\to \c[t]/(t^2)$$ 
with $x_1\mapsto a_1^{(1)}t,\ \  x_i\mapsto  a_i^{(0)}a_1^{(1)}t, (i\neq 1)$,
 it corresponds to a point 
$$(a_1^{(1)},  a_2^{(0)}a_1^{(1)},\ldots, a_N^{(0)}a_1^{(1)})$$
in $\bA^N\simeq C_{X,x}$.
Since $b_1(\alpha)$ is not a trivial jet, we have $a_1^{(1)}\neq 0$ and the  image 
of $(a_1^{(1)},  a_2^{(0)}a_1^{(1)},\ldots, a_N^{(0)}a_1^{(1)})$ by $\rho$ 
has the homogeneous coordinates $(1:a_2^{(0)}:\ldots:a_N^{(0)})$ that represents $P$.
\end{proof}

\begin{thm}
\label{main}
For a singularity $(X,x)$ the following are equivalent:
\begin{enumerate}
\item[(i)] $\hmld(x; X, \o_X)= n$;
\item[(ii)] $\lambda_m=0$ for every $m\in \n$;
\item[(iii)] $\lambda^0_m=0$ for every $m\in \n$;
\item[(iv)] $\lambda^0_1=0$;
\item[(v)] The tangent cone of $(X,x)$ has a reduced irreducible component,
where a reduced irreducible component means an irreducible component
which is reduced at the generic point.
\end{enumerate}
\end{thm}

\begin{proof} The implications (ii)$\Rightarrow$(i) and (iii)$\Rightarrow$(i) are trivial
by Lemma \ref{lambda}, (iii).
The implication (i)$\Rightarrow$(ii) follows from $0=\lambda=\lambda_m$
for $m\gg 0$.
Indeed, fix $m\gg 0$, then by \cite[Lemma 4.3]{DL2} and $\lambda_m=0$, for every $i<m$,
$$in\geq \dim \psi_i(\pi^{-1}(x))\geq mn-(m-i)n= in,$$
which means $\lambda_i=0$.
The implication (i)$\Rightarrow$(iii) follows similarly.
(iii)$\Rightarrow$(iv) is trivial.
(iv)$\Rightarrow$(v) is proved as follows:
By (iv) we have $\dim \psi_1(\pi^{-1}(x)\setminus(\sing X)_\infty)=n$.
Then, by Lemma \ref{reg},
$$n=\dim h_1\left((\pi_1^{\tilde Y})^{-1}(E'_{reg})\right)=
\dim b_1\left(g_1(\pi_1^{\tilde Y})^{-1}(E'_{reg})\right).$$
Therefore, $\dim\rho\circ\phi\left(b_1\left(g_1(\pi_1^{\tilde Y})^{-1}(E'_{reg})\right)\right)=n-1$.
By Lemma \ref{corres}, we have
$$\rho\circ\phi\left(b_1\left(g_1(\pi_1^{\tilde Y})^{-1}(E'_{reg})\right)\setminus 
\{0\}\right)\simeq
\pi_1^Y\left(g_1(\pi_1^{\tilde Y})^{-1}(E'_{reg})\right).$$
On the other hand, by  commutativity of the diagram, it follows
$$\pi_1^Y\left(g_1(\pi_1^{\tilde Y})^{-1}(E'_{reg})\right)=g(E'_{reg}).$$
Now we obtain $$\dim g(E'_{reg})=n-1.$$
This shows that $E=b^{-1}(x)$ is reduced at an irreducible component,
which implies (v).
The proof of (v)$\Rightarrow$(iii) is the following: 
  Let $C$ be the reduced irreducible component of $E=b^{-1}(x)$.
  Then, by restricting $Y$ and $C$ by appropriate open subsets, we may assume that a pair $(Y, C)$ is a non-singular pair of a variety and a divisor by \cite[Chapter 0, 17.1.8]{ega}.
  As $\psi_m((\pi^Y)^{-1}(C))=(\pi_m^Y)^{-1}(C)$ is of dimension $mn+n-1$ and
   in general $\dim \overline{\psi_m(\pi^{-1}(x)\setminus (\sing X)_\infty)}\leq mn$, 
  it is sufficient to show that the general fiber of the morphism
  $$b_m:\psi_m((\pi^Y)^{-1}(C))\to \overline{\psi_m(\pi^{-1}(x)\setminus (\sing X)_\infty)}$$
  onto its image 
  has  dimension $n-1$.
  Hence, we have only to prove that $$\dim b_m^{-1}b_m(\alpha)=n-1$$
  for a general $\alpha\in (\pi_m^Y)^{-1}(C).$

   Let $(X,x)\subset (\bA^N,0)$ and let $\tilde b: \tilde \bA^N\to \bA^N$ be the blow-up at 
  the origin $0$. As $\tilde \bA^N$ is covered by $U_{(i)}=\spec\c[x_i,\frac{x_1}{x_i},\ldots, \frac{x_N}{x_i}]$ $(i=1,\ldots, N)$, we may assume that 
 $$(Y, C)\subset (U_{(1)}, Z(x_1)).$$
 Take $\alpha\in (\pi_m^Y)^{-1}(C)$.
 Then, the $m$-jet $\alpha:\spec \c[t]/(t^{m+1})\to Y\subset U_{(1)}$ corresponding to a ring homomorphism
 $$ \c\left[x_1,\frac{x_2}{x_1},\ldots, \frac{x_N}{x_1}\right]\to\c[t]/(t^{m+1}),\ \ 
 x_1\mapsto \sum_{j=1}^{m}a_1^{(j)}t^j,\ \frac{x_i}{x_1}\mapsto \sum_{j=0}^{m}a_i^{(j)}t^j (i\neq 1)$$
 has the coordinates 
 $$(\overbrace{0,a_2^{(0)},\ldots,a_N^{(0)}}^{N},\overbrace{a_1^{(1)},\ldots,a_N^{(1)}}^{N},\ldots,\overbrace{ a_1^{(m)},\ldots,a_N^{(m)}}^{N})$$
 in $(U_{(1)})_m=\bA^{(m+1)N}$.
 The reason why $a_1^{(0)}=0$ is because $\pi_m^Y(\alpha) \in C\subset Z(x_1)$.
 
 Here, as $Y$ is non-singular and $x_1$ is a member of a regular system of  parameters of $Y$ at all points in $C$, it follows 
 $a_1^{(1)}\neq 0$ for general $\alpha\in (\pi_m^Y)^{-1}(C).$
 Indeed, we can take a regular system $x_1, z_2,\ldots, z_n$ of parameters of $\o_{Y,P}$
 for a general point $P\in C$.
 Then, we obtain an $m$-jet $\c[[x_1,z_2,\ldots, z_n]]\to \c[t]/(t^{m+1})$ which sends
 $x_1$ to $\sum_{j=0}^ma_1^{(j)}t^j$  with $a_1^{(0)}=0, a_1^{(1)}\neq 0 $.
This shows that a general $\alpha\in (\pi_m^Y)^{-1}(C)$ satisfies $a_1^{(1)}\neq 0 $.

 Then, $b_m(\alpha)$ corresponds to the ring homomorphism
$$ \c\left[x_1,\ldots, {x_N}\right]\to\c[t]/(t^{m+1}),\ \ 
 x_1\mapsto \sum_{j=1}^{m}a_1^{(j)}t^j,\ {x_i}\mapsto \sum_{j=1}^{m}\sum_{l+k=j}a_1^{(l)}a_i^{(k)}t^j (i\neq 1),$$
 because $x_i=x_1\frac{x_i}{x_1}$.
 Therefore, $b_m(\alpha)$ has the coordinates
 \begin{equation}
 \label{a}
 (\overbrace{0,\ldots,0}^{N}\mid \overbrace{a_1^{(1)},a_1^{(1)}a_2^{(0)},\ldots,a_1^{(1)}a_N^{(0)}}^{N}\mid \ldots,
 \mid \overbrace{a_1^{(m)}, \sum_{l=1}^m a_1^{(l)}a_2^{(m-l)},\ldots, \sum_{l=1}^m a_1^{(l)}a_N^{(m-l)}}^{N})
 \end{equation}
 in $(\bA^N)_m=\bA^{(m+1)N}$.
 Here, in the coordinates we put  slits to clarify each block of $N$-coordinates for the convenience.
 
 Let $\beta\in b_m^{-1}b_m(\alpha)$ and let the coordinates of $\beta$ in $(U_{(1)})_m=\bA^{(m+1)N}$ be
 $$(\overbrace{0,d_2^{(0)},\ldots,d_N^{(0)}}^{N},\overbrace{d_1^{(1)},\ldots,d_N^{(1)}}^{N},\ldots,\overbrace{ d_1^{(m)},\ldots,d_N^{(m)}}^{N}).$$
 Then, $b_m(\beta)$ has the coordinates
 \begin{equation}
 \label{d}
 (\overbrace{0,\ldots,0}^{N}\mid \overbrace{d_1^{(1)},d_1^{(1)}d_2^{(0)},\ldots,d_1^{(1)}d_N^{(0)}}^{N}\mid \ldots,
 \mid \overbrace{d_1^{(m)}, \sum_{l=1}^m d_1^{(l)}d_2^{(m-l)},\ldots, \sum_{l=1}^m d_1^{(l)}d_N^{(m-l)}}^{N}).
 \end{equation}
 Assume that $\alpha\in (\pi_m^Y)^{-1}(C)$ is general, then as $b_m(\beta)=b_m(\alpha)$,
 we compare the coordinates in (\ref{a}) and (\ref{d}).
 First, by the comparison of the first coordinates of each block in (\ref{a}) and (\ref{d}), we obtain
 $$d_1^{(j)}=a_1^{(j)} \ (j=1,\ldots, m).$$
 Next, by the comparison of the second blocks in (\ref{a}) and (\ref{d}), we have
 $$d_i^{(0)}=\frac{1}{d_1^{(1)}}a_1^{(1)}a_i^{(0)}=a_i^{(0)} \ \ (i=2,\ldots , N),$$
 since $a_1^{(1)}=d_1^{(1)}\neq 0.$
 Then, by the comparison of the third blocks in (\ref{a}) and (\ref{d}), we have 
 $$d_i^{(1)}=\frac{1}{d_1^{(1)}}(a_1^{(2)}a_i^{(0)}+a_1^{(1)}a_i^{(1)}-d_1^{(2)}d_i^{(0)})
 =\frac{1}{a_1^{(1)}}(a_1^{(1)}a_i^{(1)})=a_i^{(1)}.$$
 In the similar way, we obtain successively 
               $$d_i^{(j)}=a_i^{(j)} \ \ \  (j=1,\ldots,m-1, \ i=1,\ldots, N)$$
 and
 $$d_1^{(m)}=a_1^{(m)}.$$             
 Therefore, for a general $\alpha$, a jet $\beta\in b_m^{-1}b_m(\alpha)$ has the coordinates of the form
 $$(0,a_2^{(0)},\ldots,a_N^{(0)}\mid a_1^{(1)},\ldots,a_N^{(1)}\mid \ldots\mid 
  a_1^{(m-1)},\ldots,a_N^{(m-1)}\mid
 a_1^{(m)},c_2,c_3,\ldots,c_N)$$
 for any $c_2,\ldots, c_N$ so that $\beta$ is in $Y_m$.
 Therefore,
 $$b_m^{-1}b_m(\alpha)=(\psi_{m,m-1}^Y)^{-1}(\psi_{m,m-1}^Y(\alpha))\cap
 Z(x_1^{(m)}-a_1^{(m)}).$$
 As $Y$ is non-singular of dimension $n$, we have 
 $$(\psi_{m,m-1}^Y)^{-1}(\psi_{m,m-1}^Y(\alpha))\simeq \bA^n.$$ 
 Again as $Y$ is non-singular and $x_1$ is a member of a regular system of  parameters at a point in $C$, the equation 
   $x_1^{(m)}-a_1^{(m)}=0$ is not trivial on this space $\bA^n$.
   Indeed, we can take a regular system $x_1,z_2,\ldots, z_n$ of parameters of $\o_{Y, P}$,
   where $P=\pi_1^Y(\alpha)$.
   Let an $m$-jet $$\alpha':\widehat{\o}_{Y,P}=\c[[x_1,z_2,\ldots, z_n]]\to \c[t]/(t^{m+1})$$
   be defined by 
   $$x_1\mapsto \sum_{j=1}^{m-1}a_1^{(j)}t^j+at^m\ \ \mbox{with}\ \ a\neq a_1^{(m)},$$
   $$z_i\mapsto \sum_{j=0}^{m}a_1^{(j)}t^j, (i=2,\ldots, N).$$
   Then, $\alpha'$ gives a point in $(\psi_{m,m-1}^Y)^{-1}(\psi_{m,m-1}^Y(\alpha))\simeq \bA^n$ which does not vanish by $x_1^{(m)}-a_1^{(m)}$.

  By this, $ x_1^{(m)}-a_1^{(m)}=0$ is not trivial equation in $\bA^n$, this 
  equation gives a 
  hypersurface.
   Hence $$\dim b_m^{-1}b_m(\alpha)=n-1.$$
 \end{proof}

\begin{exmp}
\label{minimal}
Let $X$ be a hypersurface in $\bA^{n+1}$ defined by an equation $f=0$.
Then $\hmld(x,X)=n$ if and only if 
the initial term $f_*$ of $f$ has an irreducible factor of multiplicity one.
\end{exmp}

\begin{exmp}
\label{non-lower}
  In general, the map $x\mapsto
\hmld(x;X,\a) $, ($x\in X$ closed point) is not necessarily lower semi-continuous.
Let $(X,P)$ be a  singularity of a variety $X$ of dimension $n$.
Assume $\hmld(P;X,\o_X)> n$ (Such an example actually exists,
because we have only to take an  hypersurface singularity whose defining polynomial has the initial term whose factors are all with multiple powers. 
See Example \ref{minimal})
Every open neighborhood $U$ of $P$ contains a non-singular closed point $y\in X$.
As 
$$\hmld (y;X,\o_X)=\mld(y;X,\o_X)=n,$$ 
the map $x\mapsto
\hmld(x;X,\a) $, ($x\in X$ closed point) is
not lower semi-continuous.
\end{exmp}

\makeatletter \renewcommand{\@biblabel}[1]{\hfill#1.}\makeatother

\end{document}